\documentclass{article}
\usepackage{amssymb}

\input{tcilatex}
\begin{document}

\begin{center}
\textbf{\ \ Navier-Stokes problems in half space with parameters }

\textbf{Veli.B. Shakhmurov}

Department of Mechanical Engineering, Okan University, Akfirat, Tuzla 34959
Istanbul, Turkey,

E-mail: veli.sahmurov@okan.edu.tr

Khazar University, Baky Azerbaijan

\textbf{Abstract}
\end{center}

The existence, uniqueness and uniformly $L^{p}$ estimates for solutions of
the parameter dependent abstract Navier-Stokes problem on half space are
derived. In application the existence, uniqueness and uniformly $L^{\mathbf{p%
}}$ estimates for solution of the Wentzell-Robin type mixed problem for
Navier-Stokes equation is established.

\textbf{Key Word:}$\mathbb{\ }$Stokes systems, Navier-Stokes equations,
Differential equations with small parameters, Semigroups of operators,
Boundary value problems, Differential-operator equations, Maximal $L^{p}$
regularity

\begin{center}
\bigskip\ \ \textbf{MSC 2010: 35xx, 35Jxx, 35Kxx, 43Axx, 47Axx}

\bigskip \textbf{1}. \textbf{Introduction}
\end{center}

We will consider the initial boundary value problems (IBVP) for
Navier-Stokes equation (NSE) with small parameter 
\begin{equation}
\frac{\partial u}{\partial t}-\bigtriangleup _{\varepsilon }u+\left(
u.\nabla \right) u+\nabla \varphi +Au=f\left( x,t\right) ,\text{ }\func{div}%
u=0,\text{ }  \tag{1.1}
\end{equation}%
\begin{equation}
\sum\limits_{i=0}^{\nu }\varepsilon _{n}^{\sigma _{i}}\alpha _{i}\frac{%
\partial ^{i}u}{\partial x_{n}^{i}}\left( x^{\prime },0,t\right) =0\text{, }%
\nu \in \left\{ 0,1\right\} ,  \tag{1.2}
\end{equation}%
\begin{equation}
\text{ }u\left( x,0\right) =a\left( x\right) ,\text{ }x\in R_{+}^{n},\text{ }%
t\in \left( 0,T\right) ,  \tag{1.3}
\end{equation}

where%
\[
R_{+}^{n}=\left\{ x\in R^{n},\text{ }x_{n}>0,\text{ }x=\left( x^{\prime
},x_{n}\right) ,\text{ }x^{\prime }=\left( x_{1},x_{2},...,x_{n-1}\right)
\right\} ,
\]%
\[
\bigtriangleup _{\varepsilon }u=\dsum\limits_{k=1}^{n}\varepsilon _{k}\frac{%
\partial ^{2}u}{\partial x_{k}^{2}},\text{ }\sigma _{i}=\frac{1}{2}\left( i+%
\frac{1}{q}\right) \text{, }q\in \left( 1,\infty \right) \text{,}
\]%
$\alpha _{i}$ are complex numbers, $\varepsilon =\left( \varepsilon
_{1},\varepsilon _{2},...,\varepsilon _{n}\right) $, $\varepsilon _{k}$ are
small positive parameters and $A$ is a linear operator in a Banach space $E$%
. Here%
\[
u=u_{\varepsilon }\left( x,t\right) =\left( u_{1}\left( x,t\right)
,u_{2}\left( x,t\right) ,...,u_{n}\left( x,t\right) \right) ,\text{ }%
u_{k}\left( x,t\right) =u_{k\varepsilon }\left( x,t\right) 
\]%
and $\varphi =\varphi \left( x,t\right) $ are represent the $E$-valued
unknown velocity and pressure like functions, respectively; $f=\left(
f_{1}\left( x,t\right) ,f_{2}\left( x,t\right) ,...,f_{n}\left( x,t\right)
\right) $ and $a$ represent a given $E$-valued external force and the
initial velocity. In this work, we show the uniform existence and uniqueness
of the stronger local and global solution of the Navier-Stokes problem with
small parameter $\left( 1.1\right) -\left( 1.3\right) $. This problem is
characterized by presence abstract operator $A$ and\ a small parameters $%
\varepsilon _{k}$ which corresponds to the inverse of Reynolds number $Re$
very large for the Navier-Stokes equations. The regularity properties of
Navier-Stokes equations studied in e.g. $\left[ 4-6\right] $ and $\left[ 9-15%
\right] $. Navier-Stokes equations with small viscosity when the boundary is
either characteristic or non-characteristic have been well-studied see, e.g.
in $\left[ \text{9, 11, 21}\right] .$ Moreover, regularity properties of
differential operator equation (DOE) were investigated e.g. in $\left[ \text{%
1, 2, 16-20, 23}\right] .$ Here we consier Navier-Stokes operator equation
in a Banach space $E$. Since the Banach space $E$ is arbitrary and $A$ is a
possible linear operator, by chousing spaces $E$ and operators $A$ we can
obtained existence, uniqueness and $L^{p}$ estimates of solutions for
numerous class of Novier-Stokes type problems.

In this paper, firstly\ we prove that the Stokes problem%
\[
\frac{\partial u}{\partial t}-\bigtriangleup _{\varepsilon }u+Au+\nabla
\varphi =f\left( x,t\right) ,\text{ }\func{div}u=0,\text{ }x\in R_{+}^{n},%
\text{ }t\in \left( 0,T\right) , 
\]%
\begin{equation}
\sum\limits_{i=0}^{\nu }\varepsilon _{n}^{\sigma _{i}}\alpha _{i}\frac{%
\partial ^{i}u}{\partial x_{n}^{i}}\left( x^{\prime },0,t\right) =0\text{, }%
\nu \in \left\{ 0,1\right\} ,\text{ }u\left( x,0\right) =a\left( x\right) 
\tag{1.4}
\end{equation}%
has a unique solution $\left( u,\nabla \varphi \right) $ for $f\in
L^{p}\left( 0,T;L^{q}\left( R_{+}^{n};E\right) \right) =B\left( p,q\right) ,$
$p,q\in \left( 1,\infty \right) $ and the following uniform estimate holds%
\[
\left\Vert \frac{\partial u}{\partial t}\right\Vert _{B\left( p,q\right)
}+\dsum\limits_{k=1}^{n}\left\Vert \varepsilon _{k}\frac{\partial ^{2}u}{%
\partial x_{k}^{2}}\right\Vert _{B\left( p,q\right) }+\left\Vert
Au\right\Vert _{B\left( p,q\right) }+\left\Vert \nabla \varphi \right\Vert
_{B\left( p,q\right) }\leq 
\]%
\[
C\left( \left\Vert f\right\Vert _{B\left( p,q\right) }+\left\Vert
a\right\Vert _{B_{p,q}^{2-\frac{2}{p}}}\right) 
\]%
with $C=C\left( T,p,q\right) $ independent of $f$ and $\varepsilon .$

Then, by following Kato and Fujita $\left[ 6,\text{ }10\right] $ method and
using the above uniform coercive estimate for Stokes problem we derive a
local a priori estimates for solutions of $\left( 1.1\right) -\left(
1.3\right) $, i.e., we prove that for $\gamma <1$ and $\delta \geq 0$ such
that $\frac{n}{2q}-\frac{1}{2}\leq \gamma ,$ $-\gamma <\delta <1-\left\vert
\gamma \right\vert ,$ $a\in D\left( O_{\varepsilon q}^{\gamma }\right) $
there is $T_{\ast }\in \left( 0,T\right) $ independent of $\varepsilon
_{k}\in \left( 0\right. ,\left. 1\right] $ such that $\left\Vert
O_{\varepsilon q}^{-\delta }Pf\left( t\right) \right\Vert $ is continuous on 
$\left( 0,T\right) $ and satisfies $\left\Vert O_{\varepsilon q}^{-\delta
}Pf\left( t\right) \right\Vert =o\left( t^{\gamma +\delta -1}\right) $ as $%
t\rightarrow 0,$ then there is a local solution of $\left( 1.1\right)
-\left( 1.3\right) $ such that $u\in C\left( \left[ 0,T_{\ast }\right]
;D\left( O_{\varepsilon q}^{\gamma }\right) \right) $, $u\left( 0\right) =a,$
$u\in C\left( \left( 0\right. \left. T_{\ast }\right] ;D\left(
O_{\varepsilon q}^{\alpha }\right) \right) $ for some $T_{\ast
}>0,\left\Vert O_{\varepsilon q}^{\alpha }u\left( t\right) \right\Vert
=o\left( t^{\gamma -\alpha }\right) $ as $t\rightarrow 0$ for all $\alpha $
with $\gamma <\alpha <1-\delta $ uniformly in $\varepsilon .$ Moreover, the
solution of $\left( 1.1\right) -\left( 1.3\right) $ is unique if $u\in
C\left( \left( 0\right. \left. T_{\ast }\right] ;D\left( O_{\varepsilon
q}^{\beta }\right) \right) ,$ $\left\Vert O_{\varepsilon q}^{\alpha }u\left(
t\right) \right\Vert =o\left( t^{\gamma -\beta }\right) $ as $t\rightarrow 0$
for some $\beta $ with $\beta >\left\vert \gamma \right\vert \ $uniformly in 
$\varepsilon $. For sufficiently small date we show that, there is a global
solution of the problem $\left( 1.1\right) -\left( 1.3\right) $.
Particularly, we prove that there is a $\delta >0$ such that if $\left\Vert
a\right\Vert _{L^{q}\left( R_{+}^{n};E\right) }<\delta $, then there is a
global solution $u_{\varepsilon }$ of $\left( 1.1\right) -\left( 1.3\right) $
so that%
\[
t^{\left( 1-\frac{n}{q}\right) /2}u_{\varepsilon }\text{, }t^{\left( 1-\frac{%
n}{2q}\right) }\nabla u_{\epsilon }\in C\left( \left[ 0\right. ,\left.
\infty \right) ;L^{q}\left( R_{+}^{n};E\right) \right) \text{ for }n\leq
q\leq \infty .
\]%
Moreover, the following uniform estimates hold%
\[
\sup\limits_{t,\varepsilon _{k}}\left\Vert t^{\left( 1-\frac{n}{q}\right)
/2}u_{\varepsilon }\right\Vert _{L^{q}\left( R_{+}^{n};E\right) }\leq C,%
\text{ }\sup\limits_{t,\varepsilon _{k}}\left\Vert t^{\left( 1-\frac{n}{2q}%
\right) }\nabla u_{\varepsilon }\right\Vert _{L^{q}\left( R_{+}^{n};E\right)
}\leq C,\text{ }k=1,2,...,n
\]

In application we choose $E=L_{p_{1}}\left( \Omega \right) $ and $A$ to be
differential operator with generalized Wentzell-Robin boundary condition
defined by 
\[
D\left( A\right) =\left\{ u\in W_{p_{1}}^{2}\left( 0,1\right) ,\text{ }%
B_{j}u=Au\left( j\right) +\dsum\limits_{i=0}^{1}\alpha _{ji}u^{\left(
i\right) }\left( j\right) ,\text{ }j=0,1\right\} ,\text{ } 
\]%
\[
\text{ }Au=au^{\left( 2\right) }+bu^{\left( 1\right) }+cu, 
\]
in $\left( 1.1\right) -\left( 1.2\right) ,$ where $\alpha _{ji}$ are complex
numbers, $a,$ $b$, $c$ are complex-valued functions. Then, we obtain the
following Wentzell-Robin type mixed problem for Novier-Stokes equation

\begin{equation}
\frac{\partial u}{\partial t}-\bigtriangleup _{\varepsilon }u+\left(
u.\nabla \right) u+\nabla \varphi +a\frac{\partial ^{2}u}{\partial y^{2}}+b%
\frac{\partial u}{\partial y}+cu=f\left( x,y,t\right) ,\text{ }  \tag{1.5}
\end{equation}%
\[
\func{div}_{x}u=0,\text{ }u=u\left( x,y,t\right) ,\text{ }x\in R_{+}^{n}, 
\]

\begin{equation}
\sum\limits_{i=0}^{\nu }\varepsilon _{n}^{\sigma _{i}}\alpha _{i}\frac{%
\partial ^{i}u}{\partial x_{n}^{i}}\left( x^{\prime },0,y,t\right) =0\text{, 
}\nu \in \left\{ 0,1\right\} ,\text{ }x^{\prime }\in R^{n-1}\text{, }y\in
\left( 0,1\right)  \tag{1.6}
\end{equation}%
\begin{equation}
Au\left( x,j,t\right) +\dsum\limits_{i=0}^{1}\alpha _{ji}u^{\left( i\right)
}\left( x,j,t\right) =0\text{, }u\left( x,0\right) =a\left( x\right) . 
\tag{1.7}
\end{equation}%
\ \ \ 

Note that, the regularity properties of Wentzell-Robin type BVP for elliptic
equations were studied e.g. in $\left[ \text{7, 8}\right] $ and the
references therein. Here 
\[
\tilde{\Omega}=R_{+}^{n}\times \left( 0,1\right) ,\text{ }\mathbf{p=}\left(
p_{1},p\right) .
\]%
$L^{\mathbf{p}}\left( \tilde{\Omega}\right) $ denotes the space of all $%
\mathbf{p}$-summable complex-valued\ functions with mixed norm i.e., the
space of all measurable functions $f$ defined on $\tilde{\Omega}$, for which 
\[
\left\Vert f\right\Vert _{L^{\mathbf{p}}\left( \tilde{\Omega}\right)
}=\left( \int\limits_{R_{+}^{n}}\left( \int\limits_{0}^{1}\left\vert f\left(
x,y\right) \right\vert ^{p_{1}}dy\right) ^{\frac{p}{p_{1}}}dx\right) ^{\frac{%
1}{p}}<\infty .
\]

By using the above general abstract result, the existence, uniqueness and
uniformly $L^{\mathbf{p}}\left( \tilde{\Omega}\right) $ estimates for
solution of the problem $\left( 1.5\right) -\left( 1.7\right) $ is
obtained.\ 

Let $E$ be a Banach space and $L^{p}\left( \Omega ;E\right) $ denotes the
space of strongly measurable $E$-valued functions that are defined on the
measurable subset $\Omega \subset R^{n}$ with the norm

\[
\left\Vert f\right\Vert _{L^{p}}=\left\Vert f\right\Vert _{L^{p}\left(
\Omega ;E\right) }=\left( \int\limits_{\Omega }\left\Vert f\left( x\right)
\right\Vert _{E}^{p}dx\right) ^{\frac{1}{p}},1\leq p<\infty \ . 
\]

The Banach space\ $E$ is called an $UMD$-space if\ the Hilbert operator $%
\left( Hf\right) \left( x\right) =\lim\limits_{\varepsilon \rightarrow
0}\int\limits_{\left\vert x-y\right\vert >\varepsilon }\frac{f\left(
y\right) }{x-y}dy$\ is bounded in $L^{p}\left( R,E\right) ,$ $p\in \left(
1,\infty \right) $ (see. e.g. $\left[ \text{2, \S\ 4}\right] $). $UMD$
spaces include e.g. $L^{p}$, $l^{p}$ spaces and Lorentz spaces $L_{pq},$ $p$%
, $q\in \left( 1,\infty \right) $.

Let $E_{1}$ and $E_{2}$ be two Banach spaces. Let $B\left(
E_{1},E_{2}\right) $ denote the space of all bounded linear operators from $%
E_{1}$ to $E_{2}.$ For $E_{1}=E_{2}=E$ it will be denoted by $B\left(
E\right) .$

A linear operator\ $A$ is said to be positive in a Banach\ space $E$ with
bound $M>0$ if $D\left( A\right) $ is dense on $E$ and $\left\Vert \left(
A+\lambda I\right) ^{-1}\right\Vert _{B\left( E\right) }\leq M\left(
1+\left\vert \lambda \right\vert \right) ^{-1}$ for any $\lambda \in \left(
-\infty \right. ,\left. 0\right] $ where $I$ is the identity operator in $E$
(see e.g $\left[ \text{22, \S 1.15.1}\right] $).

The positive operator $A$ is said to be $R$-positive in a Banach space $E$
if the set $L_{A}=\left\{ \xi \left( A+\xi \right) ^{-1}\text{: }\xi \in
\left( -\infty \right. ,\left. 0\right] \right\} $, is $R$-bounded (see $%
\left[ \text{2, \S\ 4}\right] $).

The operator $A\left( s\right) $ is said to be positive in\ $E$ uniformly
with respect to papameter $s$ with bound $M>0$ if $D\left( A\left( s\right)
\right) $ is independent on $s$, $D\left( A\left( s\right) \right) $ is
dense in $E$ and $\left\Vert \left( A\left( s\right) +\lambda \right)
^{-1}\right\Vert \leq \frac{M}{1+\left\vert \lambda \right\vert }$ for all $%
\lambda \in S_{\psi },0\leq \psi <\pi $, where $M$ does not depend on $s$
and $\lambda .$

Assume $E_{0}$ and $E$ are two Banach spaces and $E_{0}$ is continuously and
densely embeds into $E$. Here $\Omega $ is a measurable set in $R^{n}$ and $%
m $ is a positive integer$.$ Let\ $W^{m,p}\left( \Omega ;E_{0},E\right) $
denote the space of all functions $u\in L^{p}\left( \Omega ;E_{0}\right) $
that have the generalized derivatives $\frac{\partial ^{m}u}{\partial
x_{k}^{m}}\in L^{p}\left( \Omega ;E\right) $ with the norm 
\[
\ \left\Vert u\right\Vert _{W^{m,p}\left( \Omega ;E_{0},E\right)
}=\left\Vert u\right\Vert _{L^{p}\left( \Omega ;E_{0}\right)
}+\sum\limits_{k=1}^{n}\left\Vert \frac{\partial ^{m}u}{\partial x_{k}^{m}}%
\right\Vert _{L^{p}\left( \Omega ;E\right) }<\infty . 
\]%
\ \ $\ $

\begin{center}
\textbf{2. Regularity properties of solutions for DOEs with parameters}
\end{center}

In this section, we consider the boundary value problem (BVP) for the
elliptic DOE with small parameters in half-space. We will derive the maximal
regularity properties of the following problem%
\begin{equation}
-\bigtriangleup _{\varepsilon }u+\left( A+\lambda \right) u=f\left( x\right)
,\text{ }x\in R_{+}^{n},  \tag{2.1}
\end{equation}%
\begin{equation}
\sum\limits_{i=0}^{\nu }\varepsilon _{n}^{\sigma _{i}}\alpha _{i}\frac{%
\partial ^{i}u}{\partial x_{n}^{i}}\left( x^{\prime },0,t\right) =0\text{, }
\tag{2.2}
\end{equation}%
where $A$ is a linear operator in $E$, $\alpha _{i}$ are complex numbers, $%
\varepsilon _{k}$ are positive and $\lambda $ is a complex parameters and 
\[
\bigtriangleup _{\varepsilon }u=\dsum\limits_{k=1}^{n}\varepsilon _{k}\frac{%
\partial ^{2}u}{\partial x_{k}^{2}},\text{ }\sigma _{i}=\frac{1}{2}\left( i+%
\frac{1}{q}\right) \text{, }\nu \in \left\{ 0,1\right\} . 
\]

By virtue of $\left[ \text{19,Theorem 2.2}\right] $ we have

\textbf{Theorem 2.1. }Let $E$ be a UMD\ space space and\ $A$ is an $R$%
-positive operator in $E$. Assume $m$ is a nonnegative number, $q\in \left(
1,\infty \right) ,$ $\alpha _{\nu }\neq 0,$ $0<t_{k}\leq 1,$ $k=1,2,...,n$.
Then for all $f\in W^{m,q}\left( R_{+}^{n};E\right) $, $\lambda \in S_{\psi
,\varkappa }$ and sufficiently large $\varkappa >0$ problem $\left(
2.1\right) -\left( 2.2\right) $ has a unique solution $u$ that belongs to $%
W^{2+m,q}\left( R_{+}^{n};E\left( A\right) ,E\right) $ and the following
coercive uniform estimate holds

\begin{equation}
\sum\limits_{k=1}^{n}\sum\limits_{i=0}^{m+2}\varepsilon _{k}^{\frac{i}{m+2}%
}\left\vert \lambda \right\vert ^{1-\frac{i}{m+2}}\left\Vert \frac{\partial
^{i}u}{\partial x_{k}^{i}}\right\Vert _{L^{q}\left( R_{+}^{n};E\right)
}+\left\Vert Au\right\Vert _{L^{q}\left( R_{+}^{n};E\right) }\leq
C\left\Vert f\right\Vert _{W^{m,q}\left( R_{+}^{n};E\right) }  \tag{2.3}
\end{equation}%
with $C=C\left( q,A\right) $ independent of $\varepsilon _{1},$ $\varepsilon
_{2}$,...,$\varepsilon _{n}$, $\lambda $ and $f.$

Consider the operator $Q_{\varepsilon }$ generated by problem $\left(
2.1\right) -\left( 2.2\right) $, i.e.,%
\[
D\left( Q_{\varepsilon }\right) =W^{2,q}\left( R_{+}^{n};L_{1\varepsilon
}\right) =\left\{ u\in W^{2,q}\left( R_{+}^{n}\right) ,\text{ }%
L_{1\varepsilon }u=0\right\} \text{, }
\]%
\[
Q_{\varepsilon }u=-\bigtriangleup _{\varepsilon }u+Au.
\]

From Theorem 2.1 we obtain the following

\textbf{Result 2.1}$.$ Suppose the conditions of Theorem 2.1 are satisfied.
For $\lambda \in S_{\psi ,\varkappa }$ there is a resolvent $\left(
Q_{\varepsilon }+\lambda \right) ^{-1}$ of the operator $Q_{\varepsilon }$
satisfying the following uniform estimate%
\[
\dsum\limits_{k=1}^{n}\dsum\limits_{i=0}^{2}\left\vert \lambda \right\vert
^{1-\frac{i}{2}}\varepsilon _{k}^{\frac{i}{2}}\left\Vert \frac{\partial ^{i}%
}{\partial x_{k}^{i}}\left( Q_{\varepsilon }+\lambda \right)
^{-1}\right\Vert _{B\left( L^{q}\left( R_{+}^{n};E\right) \right) }\leq C.
\]%
\textbf{\ }It is clear that the solution\ of the problem $\left( 2.1\right)
-\left( 2.2\right) $ depend on parameters $\varepsilon =\left( \varepsilon
_{1},\varepsilon _{2},...,\varepsilon _{n}\right) $, i.e. $u=u_{\varepsilon
}\left( x\right) .$ In view of the Theorem 2.1, we derive the properties of
the solutions $\left( 2.1\right) -\left( 2.2\right) .$ Particularly, by
resoning as $\left[ \text{19,Theorem 2.2}\right] $ we show the following:

\textbf{Corollary 2.1}. Let all conditions of the Theorem 2.1. hold. Then,
the solution of $\left( 2.1\right) -\left( 2.2\right) $ satisfies the
following uniform estimate%
\[
\dsum\limits_{k=1}^{n}\dsum\limits_{i=1}^{2}\varepsilon _{k}^{\frac{i}{2}%
}\left\Vert \frac{\partial ^{i}u}{\partial \varepsilon _{k}^{i}}\right\Vert
_{L^{q}\left( R_{+}^{n};E\right) }\leq \frac{C}{\left\vert \lambda
\right\vert }\left\Vert \left( Q_{\varepsilon }+\lambda \right) u\right\Vert
_{L^{q}\left( R_{+}^{n};E\right) }. 
\]

From Theorem 2.1 we obtain the following

\textbf{Result 2.2}$.$ For $\lambda \in S_{\psi ,\varkappa }$ there is a
resolvent $\left( Q_{\varepsilon }+\lambda \right) ^{-1}$ of the operator $%
Q_{\varepsilon }$ satisfying the following uniform estimate%
\begin{equation}
\dsum\limits_{k=1}^{n}\dsum\limits_{i=0}^{2}\left\vert \lambda \right\vert
^{1-\frac{i}{2}}\varepsilon _{k}^{\frac{i}{2}}\left\Vert \frac{\partial ^{i}%
}{\partial x_{k}^{i}}\left( Q_{\varepsilon }+\lambda \right)
^{-1}\right\Vert _{B\left( L^{q}\left( R_{+}^{n};E\right) \right) }\leq C. 
\tag{2.9}
\end{equation}

\begin{center}
\textbf{3. Initial-boundary value problems for Stokes system with small
parameters}
\end{center}

Consider the following BVP for the stationary Stoces equation with parameter 
\begin{equation}
-\bigtriangleup _{\varepsilon }u+Au+\nabla \varphi +\lambda u=f\left(
x\right) ,\text{ }\func{div}u=0,\text{ }x\in R_{+}^{n},  \tag{3.1}
\end{equation}%
\begin{equation}
L_{1\varepsilon }u=\sum\limits_{i=0}^{\nu }\varepsilon _{n}^{\sigma
_{i}}\alpha _{i}\frac{\partial ^{i}u}{\partial x_{n}^{i}}\left( x^{\prime
},0,t\right) =0\text{, }\nu \in \left\{ 0,1\right\} \text{. }  \tag{3.2}
\end{equation}

The function 
\[
u\in W_{\sigma }^{2,q}\left( R_{+}^{n};E\left( A\right) ,E,L_{1\varepsilon
}\right) =\left\{ u\in W^{2,q}\left( R_{+}^{n};E\left( A\right) ,E\right) ,%
\text{ }L_{1\varepsilon }u=0,\text{ }\func{div}u=0\right\} 
\]%
satisfying the equation $\left( 3.1\right) $ a.e. on $R_{+}^{n}$ is called
the stronger solution of the problem $\left( 3.1\right) -\left( 3.2\right) .$

Let $W^{s,q}\left( R_{+}^{n};E\right) $, $0<s<\infty $ be the $E-$valued
Sobolev space of order $s$ such that $W^{q,0}\left( R_{+}^{n},E\right)
=L^{q}\left( R_{+}^{n};E\right) .\ $For $q\in \left( 1,\infty \right) $ let $%
X_{q}=L_{\sigma }^{q}\left( R_{+}^{n},E\right) $ denote the closure of $%
C_{0\sigma }^{\infty }\left( R_{+}^{n};E\right) $ in $L^{p}\left(
R_{+}^{n};E\right) ,$ where 
\[
C_{0\sigma }^{\infty }\left( R_{+}^{n};E\right) =\left\{ u\in C_{0}^{\infty
}\left( R_{+}^{n};E\right) \text{, }\func{div}u=0\right\} . 
\]

By virtue of $\left[ 19\right] ,$ vector field $u\in L^{q}\left(
R_{+}^{n};E\right) $ has a Helmholtz decomposition, i.e. all $u\in
L^{q}\left( R_{+}^{n};E\right) $ can be uniquely decomposed as $%
u=u_{0}+\nabla \varphi $ with $u_{0}\in L_{\sigma }^{q}\left(
R_{+}^{n};E\right) $, $u_{0}=P_{q}u,$where $P_{q}=P$ is a projection
operator from $L^{q}\left( R_{+}^{n};E\right) $ to $L_{\sigma }^{q}\left(
R_{+}^{n};E\right) $ and $\varphi \in L_{loc}^{q}\left( R_{+}^{n};E\right) ,$
$\nabla \varphi \in L^{q}\left( R_{+}^{n};E\right) $ so that 
\[
\left\Vert \nabla \varphi \right\Vert _{q}\leq C\left\Vert u\right\Vert _{q}%
\text{, }\left\Vert \varphi \right\Vert _{L^{q}\left( G\cap B\right) }\leq
C\left\Vert u\right\Vert _{q} 
\]%
with $C$ independent of $u$, where $B$ is an open ball in $R^{n}$ and $%
\left\Vert u\right\Vert _{p}$ denotes the norm of $u$ in $L^{q}\left(
R_{+}^{n};E\right) .$

Then the problem $\left( 3.1\right) -\left( 3.2\right) $ can be reduced to
the following BVP%
\begin{equation}
-P\bigtriangleup _{\varepsilon }u+PAu+\lambda u=f\left( x\right) ,\text{ }%
x\in R_{+}^{n},  \tag{3.3}
\end{equation}%
\begin{equation}
L_{1\varepsilon }u=\sum\limits_{i=0}^{\nu }\varepsilon _{n}^{\sigma
_{i}}\alpha _{i}\frac{\partial ^{i}u}{\partial x_{n}^{i}}\left( x^{\prime
},0\right) =0\text{, }\nu \in \left\{ 0,1\right\} ,  \tag{3.4}
\end{equation}

Consider the parameter dependent Stokes operator $O_{\varepsilon }=$ $%
O_{\varepsilon ,q}$ generated by problem $\left( 3.3\right) -\left(
3.4\right) $, i.e., 
\[
D\left( O_{\varepsilon }\right) =W_{\sigma }^{2,q}\left( R_{+}^{n};E\left(
A\right) ,E,L_{1\varepsilon }\right) \text{, }O_{\varepsilon
}u=-P\bigtriangleup _{\varepsilon }u+PAu. 
\]

From the Rezult 2.2 we get that the operator $O_{\varepsilon }$ is positive
and generates a bounded holomorphic semigroup $S_{\varepsilon }\left(
t\right) =\exp \left( -O_{\varepsilon }t\right) $ for $t>0.$

In a similar way as in $\left[ 6\right] $ we show

\textbf{Proposition 3.1. }The following estimate holds 
\[
\left\Vert O_{\varepsilon }^{\alpha }S_{\varepsilon }\left( t\right)
\right\Vert \leq Ct^{-\alpha },\text{ } 
\]

uniformly in $\varepsilon =\left( \varepsilon _{1},\varepsilon
_{2},...,\varepsilon _{n}\right) $ for $\alpha \geq 0$ and $t>0.$

\textbf{Proof. }From Result 2.2 we obtain that the operator $O_{\varepsilon }
$ is uniformly positive in $L_{q}\left( R_{+}^{n};E\right) $, i.e. for $%
\lambda \in \in \left( -\infty \right. ,\left. 0\right] $ the following
uniform estimate holds%
\[
\left\Vert \left( O_{\varepsilon }+\lambda \right) ^{-1}\right\Vert \leq
M\left\vert \lambda \right\vert ^{-1},\text{ }
\]%
where the constant $M$ is independent of $\lambda $ and $\varepsilon .$
Then, by using Danford integral and operator calculus as in $\left[ 6\right] 
$ we obtain the assertion.

From $\left[ 19\right] $ we obtain the following result

\textbf{Theorem 3.1}$.$ Let $E$ be a a UMD\ space,\ $A$ an $R$-positive
operator in $E$, $q\in \left( 1,\infty \right) $ and $0<\varepsilon _{k}\leq
1$. Then for every $f\in L^{p}\left( 0,T;L^{q}\left( R_{+}^{n};E\right)
\right) =B\left( p,q\right) $ and $a\in B_{p,q}^{2-\frac{2}{p}},$ $p,q\in
\left( 1,\infty \right) $ there is a unique solution $\left( u,\nabla
\varphi \right) $ of the problem $\left( 1.9\right) $ and the following
uniform estimate holds%
\[
\left\Vert \frac{\partial u}{\partial t}\right\Vert _{B\left( p,q\right)
}+\dsum\limits_{k=1}^{n}\left\Vert \varepsilon _{k}\frac{\partial ^{2}u}{%
\partial x_{k}^{2}}\right\Vert _{B\left( p,q\right) }+\left\Vert
Au\right\Vert _{B\left( p,q\right) }+\left\Vert \nabla \varphi \right\Vert
_{B\left( p,q\right) }\leq 
\]

\begin{equation}
C\left( \left\Vert f\right\Vert _{B\left( p,q\right) }+\left\Vert
a\right\Vert _{B_{p,q}^{2-\frac{2}{p}}}\right)  \tag{3.5}
\end{equation}%
with $C=C\left( T,p,q\right) $ independent of $f$ and $\varepsilon .$

\begin{center}
\textbf{4}. \textbf{Existence and Uniqueness for Navier-Stokes equation with
parameters}
\end{center}

\bigskip In this section, we study the Navier-Stokes problem $\left(
1.1\right) -1.3$ in $X_{q}$. The problem $\left( 1.1\right) -\left(
1.3\right) $ can be expressed as%
\begin{equation}
\frac{du}{dt}+O_{\varepsilon }u=Fu+Pf,\text{ }u\left( 0\right) =0,\text{ }%
t>0,\text{ }Fu=-P\left( u,\nabla \right) u.  \tag{4.1}
\end{equation}%
We consider this equation in integral form%
\begin{equation}
u\left( t\right) =S_{\varepsilon }\left( t\right)
a+\dint\limits_{0}^{t}S_{\varepsilon }\left( t-s\right) \left[ Fu\left(
s\right) +Pf\left( s\right) \right] ds,\text{ }t>0.  \tag{4.2}
\end{equation}

For the proving the main result we need the following lemma which is
obtained from $\left[ \text{4, Theorem 2}\right] .$

\textbf{Lemma 4.1.} Let $E$ be a a UMD\ space,\ $A$ an $R$-positive operator
in $E$, $q\in \left( 1,\infty \right) $ and $0<\varepsilon _{k}\leq 1$. For
any $0\leq \alpha \leq 1$ the domain $D\left( O_{\varepsilon }^{\alpha
}\right) $ is the complex interpolation space $\left[ X_{q},D\left(
O_{\varepsilon }\right) \right] _{\alpha ,}.$

\textbf{Lemma 4.2. }Let $E$ be a a UMD\ space,\ $A$ an $R$-positive operator
in $E$, $q\in \left( 1,\infty \right) $ and $0<\varepsilon _{k}\leq 1$. For
each $k=1,2,...,n$ the operator $u\rightarrow O_{\varepsilon }^{-\frac{1}{2}%
}P\left( \frac{\partial }{\partial x_{k}}\right) u$ extends uniquely to a
uniformly bounded linear operator from $L^{q}\left( R_{+}^{n};E\right) $ to $%
X_{q}.$

\textbf{Proof. }Since $O_{\varepsilon }$ is a positive operator, it has a
fractional powers $O_{\varepsilon }^{\alpha }.$ From the Lemma 4.1 It
follows that the domain $D\left( O_{\varepsilon }^{\alpha }\right) $ is
continuously embedded in $X_{q}\cap H_{q}^{2\alpha }\left( R_{+}^{n};E\left(
A\right) ,E\right) $ for any $\alpha >0$. Then by using the duality argument
and due to uniform positivity of $O_{\varepsilon }^{\frac{1}{2}}$\ we obtain
the following uniformly in $\varepsilon $ estimate holds%
\begin{equation}
\left\Vert O_{\varepsilon }^{-\frac{1}{2}}P\left( \frac{\partial }{\partial
x_{k}}\right) u\right\Vert _{L^{q}\left( R_{+}^{n};E\right) }\leq
C\left\Vert u\right\Vert _{X_{q}}.  \tag{4.3}
\end{equation}

\textbf{\ }

By reasoning as in $\left[ 3\right] $ we obtain the following

\textbf{Lemma 4.3. }Let $E$ be a a UMD\ space,\ $A$ an $R$-positive operator
in $E$, $q\in \left( 1,\infty \right) $ and $0<\varepsilon _{k}\leq 1$. Let $%
0\leq \delta <\frac{1}{2}+\frac{n}{2}\left( 1-\frac{1}{q}\right) .$ Then the
following estimate holds%
\[
\left\Vert O_{\varepsilon }^{-\delta }P\left( u,\nabla \right) \upsilon
\right\Vert _{q}\leq M\left\Vert O_{\varepsilon }^{\theta }u\right\Vert
_{q}\left\Vert O_{\varepsilon }^{\sigma }u\right\Vert _{q} 
\]%
uniformly in $\varepsilon =\left( \varepsilon _{1},\varepsilon
_{2},...,\varepsilon _{n}\right) $ with constant $M=M\left( \delta ,\theta
,q,\sigma \right) $ provided that $\theta >0$, $\sigma >0,$ $\sigma +\delta >%
\frac{1}{2}$ and 
\[
\theta +\sigma +\delta >\frac{n}{2q}+\frac{1}{2}. 
\]

\textbf{Proof. }Assume that $0<\nu <\frac{n}{2}\left( 1-\frac{1}{q}\right) $.%
\textbf{\ }Since $D\left( O_{\varepsilon }^{\alpha }\right) $ is
continuously embedded in $X_{q}\cap H_{q}^{2\alpha }\left( R_{+}^{n};E\left(
A\right) ,E\right) $ and $L^{q^{\prime }}\left( R_{+}^{n};E\right) \cap
X_{q^{\prime }}$ is the same as $X_{s^{\prime }},$ by Sobolev imbedding
theorem we obtain that the operators 
\[
O_{\varepsilon ,q^{\prime }}^{-\nu }:X_{q^{\prime }}\rightarrow D\left(
O_{\varepsilon ,q^{\prime }}^{\nu }\right) \rightarrow X_{s^{\prime }} 
\]%
is bounded, where 
\[
\frac{1}{s^{\prime }}=\frac{1}{q^{\prime }}-\frac{2\nu }{n},\text{ }\frac{1}{%
q}+\frac{1}{q^{\prime }}=1. 
\]

By duality argument then, we get that the operator $u\rightarrow
O_{\varepsilon ,q}^{-\nu }$ is bounded from $X_{s}$ to $X_{q},$ where 
\[
\frac{1}{s}=1-\frac{1}{s^{^{\prime }}}=\frac{1}{q}+\frac{2\nu }{n}. 
\]

Consider first the case $\delta >\frac{1}{2}$. Since $P(u,\nabla )\upsilon $
is bilinear in $u,\upsilon $, it suffices to prove the estimate on a dense
subspace. Therefore assume that $u$ and $\upsilon $ are smooth. Since div $%
u=0$, we get%
\[
(u,\nabla )\upsilon =\dsum\limits_{k=1}^{n}\frac{\partial }{\partial x_{k}}%
\left( u_{k}\upsilon \right) . 
\]

Taking $\nu =\delta -\frac{1}{2}$and using the uniform boundednes of $%
O_{\varepsilon ,q}^{-\nu },$ from $X_{s}$ to $X_{q}$ and Lemma 4.2 for all $%
\varepsilon >0$ we obtain%
\[
\left\Vert O_{\varepsilon }^{-\delta }P\left( u,\nabla \right) \upsilon
\right\Vert _{q}=\left\Vert \varepsilon _{k}O_{\varepsilon ,q}^{\frac{1}{2}%
-\nu }\dsum\limits_{k=1}^{n}P\frac{\partial }{\partial x_{k}}\left(
u_{k}\upsilon \right) \right\Vert _{q}\leq \left\Vert \left\vert
u\right\vert \left\vert \upsilon \right\vert \right\Vert _{s}. 
\]

By assumption we can take $r$ and $\eta $ such that%
\[
\frac{1}{r}\geq \frac{1}{q}-\frac{2\theta }{n},\text{ }\frac{1}{\eta }\geq 
\frac{1}{q}-\frac{2\sigma }{n},\text{ }\frac{1}{r}+\frac{1}{\eta }=\frac{1}{s%
},\text{ }r>1,\text{ }\eta <\infty . 
\]

Since $D\left( O_{\varepsilon ,q}^{\alpha }\right) $ is continuously
embedded in $X_{q}\cap H_{q}^{2\alpha }\left( R_{+}^{n};E\left( A\right)
,E\right) ,$ then by Sobolev imbedding we get%
\[
\left\Vert \left\vert u\right\vert \left\vert \upsilon \right\vert
\right\Vert _{s}\leq \left\Vert u\right\Vert _{r}\left\Vert \upsilon
\right\Vert _{\eta }\leq M\left\Vert O_{\varepsilon ,q}^{\theta
}u\right\Vert _{r}\left\Vert O_{\varepsilon ,q}^{\sigma }\upsilon
\right\Vert _{\eta }, 
\]

i.e., we have the required result for $\delta >\frac{1}{2}$. In particular,
we get 
\[
\left\Vert O_{\varepsilon }^{-\frac{1}{2}}P\left( u,\nabla \right) \upsilon
\right\Vert _{q}\leq M\left\Vert O_{\varepsilon ,q}^{\theta }u\right\Vert
_{r}\left\Vert O_{\varepsilon ,q}^{\sigma }\upsilon \right\Vert _{\eta },%
\text{ }\theta +\beta \geq \frac{n}{2q},\text{ }\beta >0. 
\]

Similarly we obtain 
\[
\left\Vert P\left( u,\nabla \right) \upsilon \right\Vert _{q}\leq
C\left\Vert u\right\Vert _{r}\left\Vert \upsilon \right\Vert _{\eta }\leq
C\left\Vert O_{\varepsilon ,q}^{\theta }u\right\Vert _{r}\left\Vert
O_{\varepsilon ,q}^{\beta +\frac{1}{2}}\upsilon \right\Vert _{\eta } 
\]%
for $\frac{1}{r}+\frac{1}{\eta }=\frac{1}{q}$ and $\delta =0$. The above two
estimates show that the map $\upsilon \rightarrow P\left( u,\nabla \right)
\upsilon $ is a uniform bounded operator from $D\left( O_{\varepsilon
}^{\beta }\right) $ to $D\left( O_{\varepsilon }^{-\frac{1}{2}}\right) $ and
from $D\left( O_{\varepsilon }^{\beta +\frac{1}{2}}\right) $ to $X_{q}.$ By
using the Lemma 4.1 and the interpolation theory for $0\leq \delta \leq 
\frac{1}{2}$ we obtain 
\[
\left\Vert P\left( u,\nabla \right) \upsilon \right\Vert _{q}\leq
C\left\Vert O_{\varepsilon ,q}^{\theta }u\right\Vert _{r}\left\Vert
O_{\varepsilon ,q}^{\sigma }\upsilon \right\Vert _{\eta }. 
\]

By using Lemma 4.3 and iteration argument, by reasoning as in Fujita and
Kato $\left[ 6\right] $ we obtain the following

\textbf{Theorem 4.1.} Let $E$ be a a UMD\ space,\ $A$ an $R$-positive
operator in $E$, $q\in \left( 1,\infty \right) $ and $0<\varepsilon _{k}\leq
1$. Let $\gamma <1$ be a real number and $\delta \geq 0$ such that 
\[
\frac{n}{2q}-\frac{1}{2}\leq \gamma ,\text{ }-\gamma <\delta <1-\left\vert
\gamma \right\vert . 
\]

Suppose that $a\in D\left( O_{\varepsilon }^{\gamma }\right) $, and that $%
\left\Vert O_{\varepsilon }^{-\delta }Pf\left( t\right) \right\Vert $ is
continuous on $\left( 0,T\right) $ and satisfies%
\[
\left\Vert O_{\varepsilon }^{-\delta }Pf\left( t\right) \right\Vert =o\left(
t^{\gamma +\delta -1}\right) \text{ as }t\rightarrow 0. 
\]

Then there is $T_{\ast }\in \left( 0,T\right) $ independent of $\varepsilon $
and local solution of $\left( 4.1\right) $ such that

$u\in C\left( \left[ 0,T_{\ast }\right] ;D\left( O_{\varepsilon }^{\gamma
}\right) \right) $, $u\left( 0\right) =a,$ $u\in C\left( \left( 0\right.
\left. T_{\ast }\right] ;D\left( O_{\varepsilon }^{\alpha }\right) \right) $
for some $T_{\ast }>0,$ $\left\Vert O_{\varepsilon }^{\alpha }u\left(
t\right) \right\Vert =o\left( t^{\gamma -\alpha }\right) $ as $t\rightarrow
0 $ for all $\alpha $ with $\gamma <\alpha <1-\delta $ uniformly in $%
\varepsilon $. Moreover, the solution of $\left( 4.1\right) $ is unique if $%
u\in C\left( \left( 0\right. \left. T_{\ast }\right] ;D\left( O_{\varepsilon
}^{\beta }\right) \right) ,$ $\left\Vert O_{\varepsilon }^{\alpha }u\left(
t\right) \right\Vert =o\left( t^{\gamma -\beta }\right) $ as $t\rightarrow 0$
for some $\beta $ with $\beta >\left\vert \gamma \right\vert \ $uniformly in 
$\varepsilon =\left( \varepsilon _{1},\varepsilon _{2},...,\varepsilon
_{n}\right) $.

\textbf{Proof. }We introduce the following iteration scheme

\begin{equation}
u_{0}\left( t\right) =S_{\varepsilon }\left( t\right)
a+\dint\limits_{0}^{t}S_{\varepsilon }\left( t-s\right) Pf\left( s\right) ds,%
\text{ }  \tag{4.3}
\end{equation}

\[
u_{m+1}\left( t\right) =u_{0}\left( t\right)
+\dint\limits_{0}^{t}S_{\varepsilon }\left( t-s\right) Fu_{m}\left( s\right)
ds\text{, }m\geq 0. 
\]

By estimating the term $u_{0}\left( t\right) $ in $\left( 4.3\right) $ and
by using the Lemma 4.3 for $\gamma \leq \alpha <1-\delta $ we get 
\[
\left\Vert O_{\varepsilon }^{\alpha }u_{0}\left( t\right) \right\Vert \leq
\left\Vert O_{\varepsilon }^{\alpha }S_{\varepsilon }\left( t\right)
a\right\Vert +\dint\limits_{0}^{t}\left\Vert O_{\varepsilon }^{\alpha
+\delta }S_{\varepsilon }\left( t-s\right) \right\Vert \left\Vert
O_{\varepsilon }^{-\delta }Pf\left( s\right) \right\Vert ds\leq 
\]%
\[
\left\Vert O_{\varepsilon }^{\alpha }S_{\varepsilon }\left( t\right)
a\right\Vert +C_{\alpha +\delta }\dint\limits_{0}^{t}\left\Vert \left(
t-s\right) \right\Vert ^{-\left( \alpha +\delta \right) }\left\Vert
O_{\varepsilon }^{-\delta }Pf\left( s\right) \right\Vert ds\leq M_{\alpha
}t^{\gamma -\alpha } 
\]%
uniformly with respect to parameters $\varepsilon _{1},\varepsilon
_{2},...,\varepsilon _{n}$ with 
\[
M_{\alpha }=\sup\limits_{0<t\leq T,\text{ }\varepsilon >0}t^{\alpha -\gamma
}\left\Vert O_{\varepsilon }^{\alpha +\delta }S_{\varepsilon }\left(
t\right) a\right\Vert +C_{\alpha +\delta }NB\left( 1-\delta -\alpha ,\gamma
+\alpha \right) , 
\]%
where $N=\sup\limits_{0<t\leq T}t^{1-\gamma -\delta }\left\Vert
O_{\varepsilon }^{-\delta }Pf\left( t\right) \right\Vert $ and $B\left(
a,b\right) $ is the beta function. Here we suppose $\gamma +\delta >0.$ By
induction assume that $u_{m}\left( t\right) $ satisfies the following 
\begin{equation}
\left\Vert O_{\varepsilon }^{\alpha }u_{m}\left( t\right) \right\Vert \leq
M_{\alpha m}t^{\gamma -\alpha },\text{ }\gamma \leq \alpha <1-\delta . 
\tag{4.4}
\end{equation}

We shall estimate $O_{\varepsilon }^{\alpha }u_{m+1}\left( t\right) $ by
using $(5.2).$To estimate the term $\left\Vert O_{\varepsilon }^{-\delta
}Fu_{m}\left( s\right) \right\Vert $ we suppose 
\[
\theta +\sigma +\delta =1+\gamma ,\text{ }\gamma <\theta <1-\delta ,\text{ }%
\gamma <\sigma <1-\delta ,
\]%
\[
\text{ }\theta >0,\text{ }\sigma >0,\text{ }\delta +\sigma >\frac{1}{2},
\]%
so that the numbers $\theta ,$ $\sigma ,$ $\delta $ satisfy the assumptions
of Lemma 4.3. Using Lemma 4.3 and $\left( 4.4\right) ,$ we get 
\[
\left\Vert O_{\varepsilon }^{-\delta }Fu_{m}\left( s\right) \right\Vert \leq
CM_{\theta m}M_{\sigma m}s^{\gamma +\delta -1}.
\]

Therefore, we obtain 
\[
\left\Vert O_{\varepsilon }^{\alpha }u_{m}\left( t\right) \right\Vert \leq
M_{\alpha }t^{\gamma -\alpha }+M_{\alpha +\delta
}\dint\limits_{0}^{t}\left\Vert \left( t-s\right) \right\Vert ^{-\left(
\alpha +\delta \right) }\left\Vert O_{\varepsilon }^{-\delta }Fu_{m}\left(
s\right) \right\Vert ds 
\]

\[
\leq M_{\alpha m+1}t^{\gamma -\alpha } 
\]%
with 
\[
M_{\alpha m+1}=M_{\alpha }+M_{\alpha +\delta }MB\left( 1-\delta -\alpha
,\gamma +\delta \right) M_{\theta m}M_{\sigma m}. 
\]

We get the uniform estimate. So, the remaining part of proof is obtainedthe
same as in $\left[ \text{3},\text{Theorem 2.3}\right] $.

By reasoning as in $\left[ \text{6}\right] $ we obtain

\textbf{Lemma 4.4. }Let the parameter dependent operator $A_{\varepsilon }$
be uniform positive in a Banach space $E$ and $\alpha $ be a positive number
with $0<\alpha <1$. Then, the following uniform inequality holds%
\[
\left\Vert A_{\varepsilon }^{\alpha }\left( e^{-A_{\varepsilon }t}-I\right)
u\right\Vert _{E}\leq \frac{t^{\alpha }}{\alpha }\left\Vert A_{\varepsilon
}^{\alpha }u\right\Vert _{E} 
\]%
for all $u\in E.$

\textbf{Proposition 4.1. }Let $E$ be a space satisfying a multiplier
condition,\ $A$ an $R$-positive operator in $E$, $q\in \left( 1,\infty
\right) $ and $0<\varepsilon _{k}\leq 1$. Let $u$ be the solution given by
Theorem 4.1. Then $O_{\varepsilon }^{\alpha }u$ for $\gamma <\alpha
<1-\delta $ is uniform H\"{o}lder continuous on every interval $\left[ \eta
,T_{\ast }\right] $, $0<\eta <T_{\ast }$ for all parameters $\varepsilon
_{k}>0.$

\textbf{Proof.}\ It suffices to prove the H\"{o}lder continuity of $%
O_{\varepsilon }^{\alpha }\upsilon $, where%
\[
\upsilon \left( t\right) =\dint\limits_{0}^{t}S_{\varepsilon }\left(
t-s\right) \left[ Fu\left( s\right) +Pf\left( s\right) \right] ds. 
\]

Using the Lemma 4.4 we get the uniform estimate%
\[
\left\Vert \left( e^{-hO_{\varepsilon }}-I\right) O_{\varepsilon }^{-\alpha
}\right\Vert _{B\left( E\right) }\leq \frac{h^{\alpha }}{\alpha },\text{ }%
h>0. 
\]

Then as a similar way as in $\left[ \text{3, Proposition 2.4}\right] $ we
obtain the assertion.

\textbf{Theorem 4.2.} Let $E$ be a a UMD\ space,\ $A$ an $R$-positive
operator in $E$, $q\in \left( 1,\infty \right) $ and $0<\varepsilon _{k}\leq
1$. Assume $Pf:\left( 0\right. \left. T_{\ast }\right] \rightarrow X_{q}$ is
H\"{o}lder continuous on each subinterval $\left[ \eta ,T_{\ast }\right] .$
Then, the solution of $\left( 4.2\right) $ given by Theorem 4.1 satisfies
equation $\left( 4.1\right) $ for all parameters $\varepsilon _{k}>0.$
Moreover, $u\in D\left( O_{\varepsilon }\right) $ for $t\in \left( 0\right.
\left. T_{\ast }\right] $.

\textbf{Proof. }It suffices to show H\"{o}lder continuity of $Fu\left(
t\right) $ on each interval $\left[ \eta ,T_{\ast }\right] .$ It is clear to
see that $u\left( \eta \right) \in X_{q}$ and%
\[
u\left( t\right) =S_{\varepsilon }\left( t\right) u\left( \eta \right)
+\dint\limits_{0}^{t}S_{\varepsilon }\left( t-s\right) \left[ Fu\left(
s\right) +Pf\left( s\right) \right] ds,\text{ }t\in \left[ \eta ,T_{\ast }%
\right] . 
\]

Since $Pf$ is continuous on $\left[ \eta ,T_{\ast }\right] $ we get 
\[
\left\Vert Pf\left( t\right) \right\Vert =o\left( t-\eta \right) ^{-\alpha },%
\text{ }t\rightarrow \eta \text{, }\alpha >0. 
\]

The uniqueness of $u\left( t\right) $, ensured by Theorem 4.1, implies the
following estimates%
\[
C\left( \left[ \eta ,T_{\ast }\right] ;D\left( O_{\varepsilon }^{\nu
}\right) \right) \cap C\left( \left( \eta \right. ,\left. T_{\ast }\right]
;D\left( O_{\varepsilon }^{\alpha }\right) \right) ,
\]%
\[
O_{\varepsilon }^{\alpha }\left\Vert u\left( t\right) \right\Vert =o\left(
t-\eta \right) ^{\nu -\alpha },\text{ }t\rightarrow \eta \text{, }\nu
<\alpha <1
\]%
uniformly in $\varepsilon _{k},$ where $\nu =\max \left\{ \gamma ,0\right\} .
$ So, by Proposition 5.1, $O_{\varepsilon }^{\alpha }u\left( t\right) $ is
continuous on every subinterval $\left[ \eta ,T_{\ast }\right] .$ Since we
can choose $\theta $, $\sigma $ so that 
\[
\theta +\sigma =1+\nu \text{, }\nu <\theta <1,\text{ }\max \left\{ \gamma ,%
\frac{1}{2}\right\} <\sigma <1.\text{ }
\]

Lemma 4.2 implies that $Fu\left( t\right) $ is H\"{o}lder continuous on
every interval $\left[ \eta ,T_{\ast }\right] .$

\begin{center}
\textbf{5}. \textbf{Regularity properties }
\end{center}

The purposes of this section is to show that the solutions of the equation $%
\left( 1.1\right) $ are smooth if the data are smooth. For simplicity, we
assume $Pf=0$. The proof when $Pf\neq 0$ is the same. Consider first all of
the Stokes problem $\left( 3.3\right) -\left( 3.4\right) .$

By reason\i ng as in $\left[ \text{6, Lemma 2.14}\right] $ we obtain

\textbf{Lemma 5.1. }Let $E$ be a a UMD\ space,\ $A$ an $R$-positive operator
in $E$, $q\in \left( 1,\infty \right) $ and $0<\varepsilon _{k}\leq 1$.%
\textbf{\ }Let $f\in C^{\mu }\left( \left[ 0,T\right] ;X_{q}\right) $, for
some $\mu \in \left( 0,1\right) .$ Then for every $\eta \in \left( 0,\mu
\right) $ we have%
\[
\upsilon \left( t\right) =\dint\limits_{0}^{t}S_{\varepsilon }\left(
t-s\right) f\left( s\right) ds\ \in C^{\eta }\left( \left( 0,\right. \left.
T \right] ;D\left( O_{\varepsilon }\right) \right) \cap C^{1+\eta }\left(
\left( 0,\right. \left. T\right] ;X_{q}\right) . 
\]

In a similar way as Lemma 3.3, 3.6.,3.7 in $\left[ \text{3}\right] $ we
obtain, respectively:

\textbf{Lemma 5.2. }Let $E$ be a a UMD\ space,\ $A$ an $R$-positive operator
in $E$, $q\in \left( 1,\infty \right) $ and $0<\varepsilon _{k}\leq 1$.%
\textbf{\ } For $u,$ $\upsilon \in W^{m,q}\left( R_{+}^{n};E\left( A\right)
,E\right) ,$ $q\in \left( 1,\infty \right) $ the following hold:

(1) $Pu\in W^{m,q}\left( R_{+}^{n};E\left( A\right) ,E\right) \cap X_{q}$
and $\left\Vert Pu\right\Vert _{W^{m,q}\left( R_{+}^{n};E\right) }\leq
C_{m,q}\left\Vert u\right\Vert _{W^{m,q}\left( R_{+}^{n};E\right) };$

(2) for $m>\frac{n}{q}$ there exists a constant $C_{m,q}$ such that 
\[
\left\Vert P\left( u,\nabla \right) \upsilon \right\Vert _{W^{m,q}\left(
R_{+}^{n};E\right) }\leq C_{m,q}\left\Vert u\right\Vert _{W^{m,q}\left(
R_{+}^{n};E\right) }\left\Vert \upsilon \right\Vert _{W^{m+1,q}\left(
R_{+}^{n};E\right) }; 
\]

(3) when $q>n$ we have 
\[
\left\Vert P\left( u,\nabla \right) \upsilon \right\Vert _{L^{q}\left(
R_{+}^{n};E\right) }\leq C_{q}\left\Vert u\right\Vert _{W^{1,q}\left(
R_{+}^{n};E\right) }\left\Vert \upsilon \right\Vert _{W^{1,q}\left(
R_{+}^{n};E\right) }. 
\]

\textbf{Lemma 5.3. }Let $E$ be a UMD space,\ $A$ an $R$-positive operator in 
$E$, $q\in \left( 1,\infty \right) $ and $0<\varepsilon _{k}\leq 1$. Let $%
u=u_{\varepsilon }\left( t\right) $\ be solution of $\left( 4.2\right) $ for 
$Pf=0,$ then $u\in C^{\mu }\left( \left( 0,\right. \left. T\right] ;D\left(
O_{\varepsilon }\right) \right) $ and $\frac{du}{dt}\in C^{\mu }\left(
\left( 0,\right. \left. T\right] ;X_{q}\right) $ for $\mu \in \left( 0,\frac{%
1}{2}\right) .$ Moreover,%
\[
Fu\in C^{\mu }\left( \left( 0,\right. \left. T\right] ;W^{1,q}\left(
R_{+}^{n};E\left( A\right) ,E\right) \right) . 
\]

\textbf{Lemma 5.4. }Let $E$ be a a UMD\ space,\ $A$ an $R$-positive operator
in $E$, $q\in \left( 1,\infty \right) $ and $0<\varepsilon _{k}\leq 1$. 
\textbf{\ }Let $u=u_{\varepsilon }\left( t\right) $\ be solution of $\left(
4.2\right) $ for $Pf=0,$ then $u\in C^{\mu }\left( \left( 0,\right. \left. T%
\right] ;D\left( O_{\varepsilon }^{\frac{1}{2}}\right) \right) $ for $\mu
\in \left( 0,\frac{1}{2}\right) .$

Now by reasoning as in $\left[ \text{3, Proposition 3.5 }\right] $ we can
state the following

\bigskip \textbf{Proposition 5.1. }Let $E$ be a a UMD space,\ $A$ an $R$%
-positive operator in $E$, $q\in \left( 1,\infty \right) $ and $%
0<\varepsilon _{k}\leq 1$.\textbf{\ }Let $E$ be Banach algebra, $q>n$ and $%
a\in X_{q}.$ Suppose that the solution $u=u_{\varepsilon }\left( t\right) $
of $\left( 4.2\right) $ for $Pf=0$ given by Theorem 4.1 exists on $\left[ 0,T%
\right] .$ Then $u\in C^{\infty }\left( R_{+}^{n}\times \left[ 0,T\right]
;E\right) .$

\textbf{Proof. }The solution $u=u_{\varepsilon }\left( t\right) $ of $\left(
4.2\right) $ for $Pf=0$ given by Theorem 4.1 \ is expressed as 
\begin{equation}
u\left( t\right) =S_{\varepsilon }\left( t\right)
a+\dint\limits_{0}^{t}S_{\varepsilon }\left( t-s\right) Fu\left( s\right) ds,%
\text{ }t>0,  \tag{5.1}
\end{equation}%
where $Fu=-P\left( u,\nabla \right) u.$ From $\left( 5.1\right) $ we get 
\[
O_{\varepsilon }^{\frac{1}{2}}u\left( t\right) =S_{\varepsilon }\left(
t-\eta \right) O_{\varepsilon }^{\frac{1}{2}}u\left( \eta \right)
+\dint\limits_{\eta }^{t}O_{\varepsilon }S_{\varepsilon }\left( t-s\right)
O_{\varepsilon }^{-\frac{1}{2}}Fu\left( s\right) ds,\text{ }t>0= 
\]%
\[
S_{\varepsilon }\left( t-\eta \right) O_{\varepsilon }^{\frac{1}{2}}u\left(
\eta \right) +\upsilon \left( t\right) ,\text{ }\upsilon \left( t\right)
=\upsilon _{\varepsilon }\left( t\right) =\dint\limits_{\eta
}^{t}O_{\varepsilon }S_{\varepsilon }\left( t-s\right) O_{\varepsilon }^{-%
\frac{1}{2}}Fu\left( s\right) ds. 
\]

Since $S_{\varepsilon }\left( t-\eta \right) O_{\varepsilon }^{\frac{1}{2}%
}u\left( \eta \right) \in C^{\infty }\left( \left( \delta ,\right. \left. T%
\right] ;X_{q}\right) $ and $0<\eta <T,$ we will examining only $\upsilon
\left( t\right) $. Integrating by parts, we obtain%
\begin{equation}
\upsilon \left( t\right) =\dint\limits_{\eta }^{t}\frac{d}{ds}S_{\varepsilon
}\left( t-s\right) O_{\varepsilon }^{-\frac{1}{2}}Fu\left( s\right)
ds=\varepsilon O_{\varepsilon }^{-\frac{1}{2}}Fu\left( t\right) -  \tag{5.2}
\end{equation}%
\[
S_{\varepsilon }\left( t-\eta \right) O_{\varepsilon }^{\frac{1}{2}}Fu\left(
\delta \right) -\dint\limits_{\eta }^{t}S_{\varepsilon }\left( t-s\right)
O_{\varepsilon }^{-\frac{1}{2}}\frac{d}{ds}\left( Fu\right) \left( s\right)
ds. 
\]

Moreover, since $u\left( s\right) \in D\left( O_{\varepsilon }\right) $ for
all $\varepsilon _{k}>0$, $0<s\leq T$, we have 
\[
\left( Fu\right) \left( s\right) =-\dsum\limits_{k=1}^{n}P\left( \frac{%
\partial }{\partial x_{k}}\right) \left[ u_{k}\left( s\right) u\left(
s\right) \right] ,
\]%
where $u\left( s\right) =\left( u_{1}\left( s\right) ,u_{2}\left( s\right)
,...,u_{n}\left( s\right) \right) ,$ $u_{k}=u_{k\varepsilon }.$ Hence, by
Lemma 4.1 we get the following uniform estimate%
\[
\left\Vert O_{\varepsilon }^{-\frac{1}{2}}\frac{d}{ds}Fu\right\Vert
_{X_{q}}=\left\Vert \dsum\limits_{k=1}^{n}O_{\varepsilon }^{-\frac{1}{2}%
}P\left( \frac{\partial }{\partial x_{k}}\right) \left[ \frac{%
du_{k\varepsilon }}{ds}u_{\varepsilon }+u_{k\varepsilon }\frac{%
du_{\varepsilon }}{ds}\right] \right\Vert _{X_{q}}
\]%
\[
\leq C\left\Vert u_{\varepsilon }\right\Vert _{L^{\infty }\left(
R_{+}^{n};E\right) }\left\Vert \frac{du_{\varepsilon }}{ds}\right\Vert
_{X_{q}}\leq C\left\Vert O_{\varepsilon }^{\frac{1}{2}}u_{\varepsilon
}\right\Vert _{X_{q}}\left\Vert \frac{du_{\varepsilon }}{ds}\right\Vert
_{X_{q}}.
\]%
This estimates together with Lemma 5.3 shows that 
\[
O_{\varepsilon }^{-\frac{1}{2}}\frac{d}{ds}Fu\in C^{\mu }\left( \left(
0,\right. \left. T\right] ;X_{q}\right) .
\]%
Lemma 5.1 and Lemma 5.2 now imply that 
\[
\frac{d\upsilon }{dt}\in C^{\mu }\left( \left( 0,\right. \left. T\right]
;X_{q}\right) .
\]

Since $D\left( O_{\varepsilon }^{\frac{1}{2}}\right) \subset W^{1,q}\left(
R_{+}^{n};E\left( A\right) ,E\right) ,$ Corollary 5.1, Lemmas 5.3, 5.4 and
the identity $u\left( t\right) =O_{\varepsilon }^{\frac{1}{2}}\left( Fu-%
\frac{du}{dt}\right) $ imply 
\[
u\in C^{\mu }\left( \left( 0,\right. \left. T\right] ;W^{3,q}\left(
R_{+}^{n};E\left( A\right) ,E\right) \right) .
\]

Then the proof will be completed as in $\left[ \text{ 3, Proposition 3.5}%
\right] $ by using the induction.

Now we can state the main result of this section

\textbf{Theorem 5.1. }Let $E$ be a a UMD\ space,\ $A$ an $R$-positive
operator in $E$, $q\in \left( 1,\infty \right) $ and $0<\varepsilon _{k}\leq
1$.\textbf{\ }Let $E$ be Banach algebra and $a\in X_{q}.$ Suppose that the
solution $u=u_{\varepsilon }\left( t\right) $ of $\left( 4.2\right) $ for $%
PF=0$ given by Theorem 4.1 exists on $\left[ 0,T\right] .$ Then $u\in
C^{\infty }\left( R_{+}^{n}\times \left[ 0,T_{\ast }\right] ;E\right) .$

\textbf{Proof. }For $q>n$ the assertion is obtained from the Proposition
5.1. Let us show that the assertion is valid for $1<q\leq n.$ Indeed, the
solution $u=u_{\varepsilon }\left( t\right) $ of $\left( 5.2\right) $ for $%
PF=0$ given by Theorem 4.1 satisfies the equation $\left( 5.1\right) $ on
every subinterval $\left[ \eta ,T_{\ast }\right] ,$ $0<\eta <T$. Theorem 4.2
shows that $u_{\varepsilon }\left( \eta \right) \in D\left( O_{\varepsilon
}\right) .$ Since $0\leq \frac{n}{2q}-\frac{1}{2}\leq \gamma <1$, we have $%
D\left( O_{\varepsilon }^{\gamma }\right) \subset X_{n}$ so that $D\left(
O_{\varepsilon }\right) \subset X_{s}$ for some $s>n.$ By $\left( 4.2\right) 
$ this means that we may assume $q>n$ and $a\in X_{q}$.

\begin{center}
\textbf{6. Existence of global solutions}
\end{center}

In this section, we prove the existence and estimate of global solution of
the problem $\left( 1.1\right) -\left( 1.3\right) .$ The proofs of these
theorems are based on the theory of holomorphic semigroups and fractional
powers of generators. We assume for simplicity that $f=0$, although it is
not difficult to include nonzero $f$ under appropriate conditions. The main
result is the following

\textbf{Theorem 6.1. }Let $E$ be a UMD\ space,\ $A$ an $R$-positive operator
in $E$, $q\in \left( 1,\infty \right) $ and $0<\varepsilon _{k}\leq 1$%
\textbf{\ } and $a\in L^{q}\left( R_{+}^{n};R^{n}\right) .$ There is a $T>0$
and a unique solution $u=u_{\varepsilon }$ of $\left( 1.1\right) -\left(
1.3\right) $ so that $t^{\left( 1-\frac{n}{q}\right) /2}u\in C\left( \left[
0\right. ,\left. T\right) ;L^{q}\left( R_{+}^{n};E\right) \right) $ for $%
n\leq q\leq \infty $ and $t^{\left( 1-\frac{n}{2q}\right) }\nabla u\in
C\left( \left[ 0\right. ,\left. T\right) ;L^{q}\left( R_{+}^{n};E\right)
\right) $ for $n\leq q<\infty $. Moreover, the following estimates hold%
\[
\sup\limits_{t\in \left[ 0\right. ,\left. T\right) ,\varepsilon
_{k}>0}\left\Vert t^{\left( 1-\frac{n}{q}\right) /2}u_{\varepsilon
}\right\Vert _{L^{q}}\leq C,\text{ }\sup\limits_{t\in \left[ 0\right.
,\left. T\right) ,\varepsilon _{k}>0}\left\Vert t^{\left( 1-\frac{n}{2q}%
\right) }\nabla u_{\varepsilon }\right\Vert _{L^{q}}\leq C. 
\]

\textbf{Proof. }The solution $u=u_{\varepsilon }\left( t\right) $ of $\left(
4.2\right) $ for $Pf=0$ given by Theorem 4.1 \ is expressed as 
\begin{equation}
u\left( t\right) =u_{0}\left( t\right) +G_{\varepsilon }u\left( t\right) , 
\tag{6.1}
\end{equation}%
where, 
\[
u_{0}\left( t\right) =S_{\varepsilon }\left( t\right) a,\text{ }%
G_{\varepsilon }u\left( t\right) =\dint\limits_{0}^{t}S_{\varepsilon }\left(
t-s\right) Fu\left( s\right) ds,\text{ }t>0. 
\]%
By applying the generalized Minkovskii inequality and by Proposition 3.1 we
can see that 
\[
\left\Vert S_{\varepsilon }\left( t\right) u\right\Vert _{L^{p}}\leq
C\varepsilon _{k}^{\frac{n}{2}\left( 1+\frac{1}{p}\right) }t^{-\frac{n}{2}%
\left( 1-\frac{1}{p}\right) }\left\Vert u\right\Vert _{L^{p}}\text{, }%
k=1,2,...,n. 
\]%
By using the above estimate we get%
\begin{equation}
\left\Vert S_{\varepsilon }\left( t\right) u\right\Vert _{L^{q}}\leq
C\varepsilon _{k}^{\frac{n}{2}\left( 2+\frac{1}{q}-\frac{1}{p}\right) }t^{-%
\frac{n}{2}\left( \frac{1}{p}-\frac{1}{q}\right) }\left\Vert u\right\Vert
_{L^{p}}\text{, }  \tag{6.2}
\end{equation}%
\begin{equation}
\left\Vert \nabla S_{\varepsilon }\left( t\right) u\right\Vert _{L^{q}}\leq
C\varepsilon _{k}^{\frac{n}{2}\left( \frac{3}{2}-\left( \frac{1}{p}-\frac{1}{%
q}\right) \right) }t^{-\left( 1+\frac{n}{2}\left( \frac{1}{p}-\frac{1}{q}%
\right) \right) }\left\Vert u\right\Vert _{L^{p}}\text{ for }1<p\leq
q<\infty .  \tag{6.3}
\end{equation}

Moreover, by using $\left( 6.1\right) $, $\left( 6.2\right) $ and \ by
applying the H\"{o}lder inequality, we get 
\begin{equation}
\left\Vert F\left( u,\upsilon \right) \right\Vert _{L^{q}}\leq C\left\Vert
u\right\Vert _{L^{r}}\left\Vert \nabla \upsilon \right\Vert _{L^{s}}\text{, }%
\frac{1}{q}=\frac{1}{r}+\frac{1}{s}.  \tag{6.4}
\end{equation}

Then \i n view of $\left( 6.1\right) $-$\left( 6.4\right) $ we obtain the
following uniform estimate%
\begin{equation}
\left\Vert G_{\varepsilon }u\right\Vert _{L^{m/\gamma }}\leq
C\dint\limits_{0}^{t}\left( t-s\right) ^{-\left( \alpha +\beta -\gamma
\right) /2}\left\Vert u\left( s\right) \right\Vert _{m/\alpha }\left\Vert
\nabla u\left( s\right) \right\Vert _{m/\beta }ds,  \tag{6.5}
\end{equation}%
\begin{equation}
\left\Vert \nabla G_{\varepsilon }u\right\Vert _{L^{m/\gamma }}\leq
C\dint\limits_{0}^{t}\left( t-s\right) ^{-\left( 1+\alpha +\beta -\gamma
\right) /2}\left\Vert u\left( s\right) \right\Vert _{m/\alpha }\left\Vert
\nabla u\left( s\right) \right\Vert _{m/\beta }ds,  \tag{6.6}
\end{equation}%
where 
\[
\alpha ,\text{ }\beta ,\text{ }\gamma >0,\text{ }\gamma \leq \alpha +\beta
<n. 
\]

Then solving the equation $\left( 6.1\right) $ by successive approximation,
starting with $u_{0}=S_{\varepsilon }\left( t\right) a$ we get 
\begin{equation}
u_{k+1}=u_{0}+G_{\varepsilon }u_{k}\text{, }u_{k}=u_{k\varepsilon }\left(
t\right) \text{, }k=0\text{, }1,2,....,\text{ }  \tag{6.7}
\end{equation}%
First by reasoning as in $\left[ \text{22, Theorem 1}\right] $ and by using $%
\left( 6.3\right) $-$\left( 6.5\right) $ we show by induction that $%
u_{k}=u_{\varepsilon k}$ exists, moreover,%
\[
t^{\left( 1-\delta \right) /2}u_{\varepsilon k}\in C\left( \left[ 0\right.
,\left. \infty \right) ;L^{n/\delta }\left( R_{+}^{n};E\right) \right) ,%
\text{ }t^{1/2}\nabla u_{\varepsilon k}\in C\left( \left[ 0\right. ,\left.
\infty \right) ;L^{n}\left( R_{+}^{n};E\right) \right) 
\]%
and for $\delta \in \left( 0,1\right) $ the following uniform estimates hold%
\begin{equation}
\sup\limits_{t,\varepsilon _{k}}\left\Vert t^{\left( 1-\delta \right)
/2}u_{\varepsilon k}\right\Vert _{L^{^{n/\delta }}}\leq M_{k},\text{ }%
\sup\limits_{t,\varepsilon _{k}}\left\Vert t^{1/2}\nabla u_{\varepsilon
k}\right\Vert _{L^{q}}\leq M_{k}^{^{\prime }}.  \tag{6.8}
\end{equation}

\ By applying $\left( 6.3\right) $-$\left( 6.5\right) $ for $q=n$ and $p=%
\frac{n}{\delta }$ we have 
\begin{equation}
M_{0}=M_{0}^{^{\prime }}=C\left\Vert a\right\Vert _{L^{n}\left(
R_{+}^{n};E\right) },  \tag{6.9}
\end{equation}%
where $C$ is a positive constant. From $\left( 6.5\right) $ and $\left(
6.7\right) $ for $n\leq p<\infty $ we obtain 
\[
\left\Vert u_{\varepsilon k+1}\right\Vert _{L^{^{p}}}\leq \left\Vert
u_{\varepsilon 0}\right\Vert _{L^{^{p}}}\leq 
\]%
\[
CM_{k}M_{k}^{^{\prime }}\dint\limits_{0}^{t}\left( t-s\right) ^{-\left(
1+\delta -n/q\right) /2}s^{-\left( 1-\delta /2\right) }ds\leq Mt^{-\left(
1-n/q\right) /2}. 
\]

It follows that $u_{\varepsilon k}\left( t\right) $ converges to a limit
function $u_{\varepsilon }$ uniformly with respect to $\varepsilon =\left(
\varepsilon _{1},\varepsilon _{2},...,\varepsilon _{n}\right) $, moreover, $%
u_{\varepsilon }\in C\left( \left[ 0\right. ,\left. T\right) ;L^{n}\left(
R_{+}^{n};E\right) \right) $ for $p=n$ and $u_{\varepsilon }$ satisfies $%
\left( 6.1\right) $ for $n<p<\infty .$

\textbf{Theorem 6.2. }Let $E$ be a a UMD\ space,\ $A$ an $R$-positive
operator in $E$, $q\in \left( 1,\infty \right) $ and $0<\varepsilon _{k}\leq
1$.\textbf{\ } There is a $\mu >0$ such that if $\left\Vert a\right\Vert
_{L^{q}\left( R_{+}^{n};E\right) }<\mu $, then there is a global solution $%
u_{\varepsilon }$ of the problem $\left( 1.1\right) -\left( 1.3\right) $, so
that $t^{\left( 1-\frac{n}{q}\right) /2}u_{\varepsilon }\in C\left( \left[
0\right. ,\left. \infty \right) ;L^{q}\left( R_{+}^{n};E\right) \right) $
for $n\leq q\leq \infty ,$ $t^{\left( 1-\frac{n}{q}\right) /2}$ and $%
t^{\left( 1-\frac{n}{2q}\right) }\nabla u_{\epsilon }\in C\left( \left[
0\right. ,\left. \infty \right) ;L^{q}\left( R_{+}^{n};E\right) \right) $
for $n\leq q<\infty $. Moreover, the following uniform estimates hold%
\begin{equation}
\sup\limits_{t,\varepsilon _{k}}\left\Vert t^{\left( 1-\frac{n}{q}\right)
/2}u_{\varepsilon }\right\Vert _{L^{q}}\leq C,\text{ }\sup\limits_{t,%
\varepsilon _{k}}\left\Vert t^{\left( 1-\frac{n}{2q}\right) }\nabla
u_{\varepsilon }\right\Vert _{L^{q}}\leq C.  \tag{6.10}
\end{equation}

\textbf{Proof. }It is clear to see from proof of Theorem 6.1 that $M_{k}$
and $M_{k}^{^{\prime }}$ are bounded by a constant $M$ if $M_{0}\leq \lambda 
$. By $\left( 7.9\right) $ this is true if $\left\Vert a\right\Vert
_{L^{q}\left( R_{+}^{n};E\right) }$ is sufficiently small. In this case, as
in $\left[ \text{10}\right] $ we prove that the sequences $t^{\left(
1-\delta \right) /2}u_{\varepsilon k}$, $t^{1/2}\nabla u_{\varepsilon k}$
are bounded on $\left( 0,\infty \right) $ uniformly in $k$ and $\varepsilon
_{1},\varepsilon _{2},...,\varepsilon _{n}$ i.e., 
\begin{equation}
\sup\limits_{t,\varepsilon _{k}}\left\Vert t^{\left( 1-\delta \right)
/2}u_{\varepsilon k}\right\Vert _{L^{^{n/\delta }}}\leq M_{1},\text{ }%
\sup\limits_{t,\varepsilon _{k}0}\left\Vert t^{1/2}\nabla u_{\varepsilon
k}\right\Vert _{L^{q}}\leq M_{2}.  \tag{6.11}
\end{equation}%
Then $\left( 6.11\right) $ is obtained from $\left( 6.10\right) .$

\textbf{Remake 6.1. }Let $E$ be a a UMD\ space,\ $A$ an $R$-positive
operator in $E$, $q\in \left( 1,\infty \right) $ and $0<\varepsilon _{k}\leq
1$.\textbf{\ }Theorem 6.2 shows that all $L^{p}$ norms of $u_{\varepsilon
}\left( t\right) $ decay as $t\rightarrow \infty $ for $p>q$ uniformly in $%
\varepsilon =\left( \varepsilon _{1},\varepsilon _{2},...,\varepsilon
_{n}\right) .$

For $p=q$ we obtain the following result

\bigskip \textbf{Theorem 6.3. }Let all conditions of Theorem 6.2 hold. Then $%
\left\Vert u_{\varepsilon }\left( t\right) \right\Vert _{p}\rightarrow 0$
uniformly in $\varepsilon $ as $t\rightarrow \infty $. More precisely, we
have%
\[
\left\Vert u_{\varepsilon }\left( t\right) -u_{0\varepsilon }\left( t\right)
\right\Vert _{p}=O\left( t^{-\frac{\delta }{2}}\right) \text{ as }%
t\rightarrow \infty , 
\]%
where, $u_{0\varepsilon }\left( t\right) =S_{\varepsilon }\left( t\right) a$
and $\delta <\min \left\{ 1,n-\frac{n}{q},\frac{n}{q}-1\right\} .$

\begin{center}
\textbf{7.} \textbf{The Wentzell-Robin type mixed problem for Novier-Stokes
equations}
\end{center}

Consider the problem $\left( 1.5\right) -\left( 1.7\right) $. Here, $W^{2,%
\mathbf{p}}\left( \tilde{\Omega}\right) $ denotes the Sobolev space with
corresponding mixed norm

\bigskip The main aim of this section is to prove the following result:

\bigskip \textbf{Theorem 7.1. }Let $a\in W^{1,\infty }\left( 0,1\right) $, $%
a\left( x\right) \geq \delta >0,$ $b,c\in L^{\infty }\left( 0,1\right) $.
Suppose the condition 7.1 hold. Let $\gamma <1$ be a real number and $\delta
\geq 0$ such that 
\[
\frac{n}{2q}-\frac{1}{2}\leq \gamma ,\text{ }-\gamma <\delta <1-\left\vert
\gamma \right\vert . 
\]

Suppose $a\in D\left( O_{\varepsilon }^{\gamma }\right) $ such that $%
\left\Vert O_{\varepsilon }^{-\delta }Pf\left( t\right) \right\Vert $ is
continuous on $\left( 0,T\right) $ and satisfies%
\[
\left\Vert O_{\varepsilon }^{-\delta }Pf\left( t\right) \right\Vert =o\left(
t^{\gamma +\delta -1}\right) \text{ as }t\rightarrow 0. 
\]

Then there is $T_{\ast }\in \left( 0,T\right) $ independent of $\varepsilon $
and local solution of $\left( 4.1\right) $ such that

$u\in C\left( \left[ 0,T_{\ast }\right] ;\right) $, $u\left( 0\right) =a,$ $%
u\in C\left( \left( 0\right. \left. T_{\ast }\right] ;D\left( O_{\varepsilon
}^{\alpha }\right) \right) $ for some $T_{\ast }>0,$ $\left\Vert
O_{\varepsilon }^{\alpha }u\left( t\right) \right\Vert =o\left( t^{\gamma
-\alpha }\right) $ as $t\rightarrow 0$ for all $\alpha $ with $\gamma
<\alpha <1-\delta $ uniformly with respect to $\varepsilon .$ Moreover, the
solution of $\left( 4.1\right) $ is unique if $u\in C\left( \left( 0\right.
\left. T_{\ast }\right] ;D\left( O_{\varepsilon }^{\beta }\right) \right) $, 
$\left\Vert O_{\varepsilon }^{\alpha }u\left( t\right) \right\Vert =o\left(
t^{\gamma -\beta }\right) $ as $t\rightarrow 0$ for some $\beta $ with $%
\beta >\left\vert \gamma \right\vert \ $uniformly in $\varepsilon =\left(
\varepsilon _{1},\varepsilon _{2},...,\varepsilon _{n}\right) $.

Then problem $\left( 1.5\right) -\left( 1.7\right) $ has a unique local
strange solution $u\in C^{\left( 2\right) }\left( \left[ 0\right. ,\left.
T_{0}\right) ;Y_{\infty }^{2,p}\right) $, where $T_{0}$ is a maximal time
interval that is appropriately small relative to $M$. Moreover, if

\[
\sup_{t\in \left[ 0\right. ,\left. T_{0}\right) }\left( \left\Vert
u\right\Vert _{Y^{2,p}}+\left\Vert u\right\Vert _{X_{\infty }}+\left\Vert
u_{t}\right\Vert _{Y^{2,p}}+\left\Vert u_{t}\right\Vert _{X_{\infty
}}\right) <\infty 
\]%
then $T_{0}=\infty .$

\ \textbf{Proof.} Let $E=L^{p_{1}}\left( 0,1\right) $. It is known $\left[ 2%
\right] $\ that $L^{p_{1}}\left( 0,1\right) $ is an $UMD$ space for $%
p_{1}\in \left( 1,\infty \right) .$ Consider the operator $A$ defined by 
\[
D\left( A\right) =W^{2,p_{1}}\left( \Omega ;B_{j}u=0\right) ,\text{ }Au=a%
\frac{\partial ^{2}u}{\partial y^{2}}+b\frac{\partial u}{\partial y}+cu. 
\]

Therefore, the problem $\left( 1.7\right) -\left( 1.8\right) $ can be
rewritten in the form of $\left( 1.1\right) -\left( 1.3\right) $, where $%
u\left( x\right) =u\left( x,.\right) ,$ $f\left( x\right) =f\left(
x,.\right) $\ are functions with values in $E=L^{p_{1}}\left( 0,1\right) .$
From $\left[ \text{7, 8}\right] $ we get that the operator $A$ generates
analytic semigroup in $L^{p_{1}}\left( 0,1\right) .$ Moreover, we obtain
that the operator $A$ is $R$-positive in $L^{p_{1}}.$ Then from Theorem 4.1
we obtain the assertion.

\bigskip

\textbf{References}

\begin{enumerate}
\item H. Amann, Linear and quasi-linear equations,1, Birkhauser, Basel 1995.

\item R. Denk, M. Hieber, J. Pruss, $R$-boundedness, Fourier multipliers and
problems of elliptic and parabolic type, Mem. Amer. Math. Soc. 2005,166
(788).

\item Y. Giga, T. Miyakava, Solutions in $L_{r}$ of the Navier-Stokes
initial value problem, Arch. Ration. Mech. Anal., 1985, 89, 267-281.

\item Y. Giga, Domains of fractional powers of the Stokes operator in $L_{r}$
spaces, Arch. Ration. Mech. Anal., 1985, 89, 251-265.

\item D. Fujiwara and H. Morimoto, An $L_{r}$-theorem of the Helmholtz
decomposition of vector fields, J. Fac. Sci. Univ. Tokyo, Sec. 1977, (I) 24,
685-700.

\item H. Fujita and T. Kato, On the Navier-Stokes initial value problem I.,
Arch. Rational Mech. Anal., 1964, 16, 269-315.

\item A. Favini, G. R. Goldstein, Jerome A. Goldstein and Silvia Romanelli,
Degenerate Second Order Differential Operators Generating Analytic
Semigroups in $L_{p}$ and $W^{1,p}$, Math. Nachr. 238 (2002), 78 -- 102.

\item V. Keyantuo, M. Warma, The wave equation with Wentzell--Robin boundary
conditions on Lp-spaces, J. Differential Equations 229 (2006) 680--697.

\item M. Hamouda and R. Temam, Some singular perturbation problems related
to the Navier-Stokes equations, Advances in deterministic and stochastic
analysis, 197-227, World. Sci. Publ., Hackensack, NJ, 2007.

\item T. Kato and H. Fujita, On the nonstationary Navier-Stokes system.
Rend. Sem. Mat. Univ. Padova, 1962, 32, 243-260.

\item N. Masmoudi, Examples of singular limits in hydrodynamics. In
Evolutionary equations. Vol. III, Handb. Differ. Equ.
Elsevier/North-Holland, Amsterdam, 2007, 195--276.

\item O. A. Ladyzhenskaya, The mathematical theory of viscous incompressible
flow. Gordon and Breach, New York, 1966

\item J-L. Lions, Mathematical topics in fluid mechanics. V. 1. The
Clarendon Press Oxford University Press, New \ York, 1996. Incompressible
models, Oxford Science Publications.

\item V. Solonnikov, Estimates for solutions of nonstationary Navier-Stokes
equations, J. Sov. Math., 1977(8), 467-529.

\item P. E. Sobolevskii, Study of Navier-Stokes equations by the methods of
the theory of parabolic equations in Banach spaces, Soviet Math. Dokl.
1964(5), 720-723.

\item V. B. Shakhmurov, The Cauchy problem for generalized abstract
Boussinesq equations, Dynamic systems and applications, 25, (2016),109-122.

\item V. B. Shakhmurov, Linear and nonlinear abstract equations with
parameters, Nonlinear Anal-Theor., 2010, v. 73, 2383-2397.

\item V. B. Shakhmurov, A. Shahmurova, Nonlinear abstract boundary value
problems atmospheric dispersion of pollutants, Nonlinear Anal-Real., 2010
v.11 (2), 932-951.

\item V. B. Shakhmurov, Stokes equations with small parameters in half
plane, Computers and Mathematics with Applications, 2014, 67, 91-115.

\item V. B. Shakhmurov, Nonlocal Novier-Stokes problem with small parameter,
Boundary Value Problems, 2013, v. 2013-107, 1-21.

\item R. Temam, Navier-Stokes Equations, North-Holland, Amsterdam, 1984.

\item H. Triebel, Interpolation theory. Function spaces. Differential
operators, North-Holland, Amsterdam, 1978.

\item S. Yakubov and Ya. Yakubov, Differential-operator Equations. Ordinary
and Partial \ Differential Equations , Chapman and Hall /CRC, Boca Raton,
2000.
\end{enumerate}

\bigskip

\end{document}